\documentclass[11pt]{amsart}
\usepackage[cp1250]{inputenc}
\usepackage[T1]{fontenc}
\usepackage[T1]{fontenc}
\usepackage[all]{xy}
\usepackage{t1enc}
\usepackage{amsfonts,amsmath,amssymb,amsthm,color,amsbsy}
\usepackage{color}
\usepackage{url}
\usepackage{bbm}
\usepackage{eucal}
\usepackage{enumerate}
\usepackage{graphicx}
\usepackage{pstricks-add}
\usepackage{pst-plot}
\usepackage{float}
\usepackage{pst-node}
\usepackage{url}
\usepackage{paralist}
\usepackage{mathrsfs}
\usepackage{colortbl}
\usepackage{bm}

\newtheorem{thm}{Theorem}[section]
\newtheorem{prop}{Proposition}[section]
\newtheorem{lem}[thm]{Lemma}
\newtheorem{cor}[thm]{Corollary}
\theoremstyle{definition}

\newtheorem{rem}[thm]{Remark}
\newtheorem{ex}[thm]{Example}
\renewcommand{\geq}{\geqslant}
\renewcommand{\leq}{\leqslant}
\renewcommand{\div}{\mathrm{div}}

\def\dist{\mathrm{dist}}
\def\la{\langle}
\def\ra{\rangle}

\def\vp{\varphi}

\def\Gr{\mathrm{Gr}}
\def\m{\multimap}
\def\f{\varphi}
\def\E{\mathbb{E}}

\def\supp{\mathrm{supp}\,}
\def\eps{\varepsilon}
\def\e{\varepsilon}
\def\pr{\right )}
\def\le{\left (}

\def\cl{\mathrm{cl}\,}
\def\I{\mathrm{I}}

\def\multi{\multimap}
\def\Q{\mathbb{Q}}
\def\R{\mathbb{R}}
\def\N{\mathbb{N}}
\def\mm{\kern +2pt \raisebox{+0.5 pt}{$\shortmid$}\kern -2pt\hbox{$\multimap$}\kern +2pt}
\def\H1{H^1\le\Omega,\mathbb{R}^n\pr}
\def\L2{L^2\le\Omega,\mathbb{R}^n\pr}
\def\h10{H^1_0\le\Omega,\mathbb{R}^n\pr}

   \newcommand{\be}{\begin{equation}}
  \newcommand{\ee}{\end{equation}}
\numberwithin{equation}{section}
\author{Wojciech Kryszewski}\thanks{The second author was partially supported by the Polish National Science Center
under grant 2013/09/B/ST1/01963.}
\address{Faculty of Mathematics and Computer Sciences, Nicolaus Copernicus University, Chopina 12/18, 87-100 Toru\'n, Poland}
\email{wkrysz@mat.umk.pl}

\author{Jakub Siemianowski}
\address{Faculty of Mathematics and Computer Sciences, Nicolaus Copernicus University, Chopina 12/18, 87-100 Toru\'n, Poland}
\email{jsiem@mat.umk.pl}

\title[Bolzano theorem and PDE]{The Bolzano mean-value theorem and partial differential equations}

\date{November 03, 2016}

\subjclass[2010]{47H10, 35J47, 47B44, 35K57}

\keywords {Steady state solution, equilibrium, fixed point, tangent cone, tangency conditions, reaction-diffusion equation, advection}

\begin{document}

\baselineskip15pt

\maketitle
\begin{abstract} We study the existence of solutions to abstract equations of the form $0=Au+F(u)$, $u\in K\subset E$, where $A$ is an abstract differential operator acting in a Banach space $E$, $K$ is a closed convex set of constraints being invariant with respect to resolvents of $A$ and perturbations are subject to different tangency condition. Such problems are closely related to the so-called Poicar\'e-Miranda theorem, being the multi-dimensional counterpart of the celebrated Bolzano intermediate value theorem. In fact our main results can and should be regarded as infinite-dimensional variants of Bolzano and Miranda-Poincar\'e theorems. Along with single-valued problems we deal with set-valued ones, yielding the existence of the so-called constrained equilibria of set-valued maps.  The abstract results are applied to show existence of (strong) steady state solutions to some weakly coupled systems of drift reaction-diffusion equations or differential inclusions of this type. In particular we get the existence of strong solutions to the Dirichlet, Neumann and periodic boundary problems for  elliptic partial differential inclusions under the presence of state constraints of different type. Certain aspects of the Bernstein theory for bvp for second order ODE are studied, too. No assumptions concerning structural coupling (monotonicity, cooperativity) are undertaken.
\end{abstract}

\section{Introduction}

The purpose of this paper is twofold. On one hand, being motivated by some concrete applications (see subsection \ref{equations}), we want to establish appropriate topological tools to study the existence of solutions to systems of $N$ partial differential equations or $N$-dimensional partial differential inclusions subject to various boundary conditions and under state constraints. The presence of such constraints is justified and explained below. This is closely related to the method of the so-called `moving rectangles' (see e.g. \cite{Smoller}) and corresponding techniques used for the study of long time behavior of evolution systems. We, however, leave aside questions concerning existence, stability and invariance of solutions of parabolic evolution equations, but in this paper we confine ourselves to elliptic equations rather and their solutions, i.e. steady state or stationary solutions to  related evolutions problems. Nevertheless the `evolution' origin of the studied steady state problems is of great importance.\\
 \indent On the other hand the proposed topological methods are closely related to problems of the existence of constrained equilibria or fixed points of abstract single- or set-valued maps, having origins in the Bolzano mean-value theorem (see subsection \ref{zeros}). This celebrated result is perhaps the most important topological device when studying one-dimensional equations of the form $f(x)=0$. This fact was extensively used and generalized by numerous authors for almost 150 years (see \cite{Mawhin 1} and \cite{Browder 1}) and various important results were established.  One of the best known statements in this direction is the Poincar\'e-Miranda theorem, which is a direct $N$-dimensional version of the Bolzano theorem. We develop the infinite-dimensional counterparts of Poincar\'e-Miranda theorem, show their relation with different branches of research concerning e.g. viability theory for differential inclusions and, finally, apply them in the context of constrained PDE.\\
 \indent The notation used throughout the paper is standard. In particular $x\cdot y$ is the scalar product of $x,y\in\R^N$ and $|x|=\sqrt{x\cdot x}$ stands for the norm of $x$. The use of function spaces ($L^p$, Sobolev etc.), linear (unbounded in general) operators in Banach spaces, $C_0$ semigroups is standard. In the paper, for the sake of generality,  we deal mostly with set-valued maps (the terminology in set-valued analysis taken after \cite{Aubin}: the symbol $\multi$ denotes a set-valued map with at least closed values). It is however important to observe that results we propose are, to the best of our knowledge, new in the single-valued case, too.\\
 \indent The paper is organized as follows: in Section \ref{second} we discuss origins of problems and motivations of main assumptions; in Section \ref{third} we establish main abstract results, while Section \ref{fourth} is devoted to applications. Section \ref{third} concludes with subsection \ref{invariance} and a  discussion of invariance issues playing an important role in the paper.

\section{The motivation}\label{second}
\subsection{Drift reaction-diffusion equations}\label{equations}

When dealing with an evolving in time multicomponent active continuous substance, whose components interact via certain reaction mechanism, such as e.g. predator-prey, activator-inhibitor, competition, reaction kinetics etc., and they all diffuse with different (in general) diffusive constants and are subject to drift or advection, i.e. a passive transfer caused by, for instance, the moving ambient media, such as gas or fluid, then the adequate model is provided by the so-called systems of drift reaction-diffusion equations (see e.g. \cite{Mih}). Such systems in general are of the form
\be\label{problem 1}
\partial_tu_i=\div(d_i\nabla u_i)+f_i(t,x,u,\nabla u),\;\; i=1,...,N\ee
or, shortly,
$$u_t=Lu+f(x,u,\nabla u)\;\;\text{where}\;\; Lu=(v_1,..,v_N), \;v_i=\div(d_i\nabla u_i)\;\;\text{for}\;\; i=1,...,N$$
along with initial and the (Neumann or Dirichlet) boundary conditions on $\partial\Omega$. Here the unknown $u=(u_1,...,u_N)$ depends on spatial variables $x=(x_1,...,x_M)\in\Omega\subset\R^M$, $\Omega$ is an open smooth domain, and time $t\in [0,T]$, $T>0$, $\nabla u=\left[\frac{\partial u_i}{\partial x_j}\right]$ is the derivative of $u$. The diffusive coefficients $d_i\in C^1(\bar\Omega)$, $i=1,...,N$, and $f=(f_1,...,f_N)$ is the source/advection term depending on $(t,x,u,\nabla u)$. Sometimes it is convenient to distinguish the advection term
\be\label{advection}f_i(t,x,u,\nabla u)=g_i(t,x,u)-\gamma^i\cdot \nabla u_i,\;\; i=1,...,N,\ee
or, shortly,
$$f(t,x,u,\nabla u)=g(t,x,u)-\Gamma u\;\;\text{where}\;\; \Gamma u=(v_1,..,v_N), \;v_i=\gamma^i\cdot \nabla u_i\;\;\text{for}\;\; i=1,...,N,$$
where drift vectors  $\gamma^i:=[\gamma^i_1,...,\gamma^i_M]\in L^\infty(\Omega,\R^M)$,  $i=1,...,N$; the second summand is responsible for the drift of the system.\\
\indent Our interest is mainly focused on ecological or chemical systems, where $u_i(x,t)$ is the concentration at $x\in\Omega$ and time $t\in [0,T]$ of the $i$-th reactant, $i=1,...,N$, contained in a bounded stirred up vessel (or reactor) $\Omega$. Clearly the initial state $u(\cdot,0)\geq 0$ on $\Omega$ and the natural expectation is that $u_i(x,t)\geq 0$ since the concentration cannot be negative. On the other hand there is a threshold value $R_i>0$  beyond which the $i$-th component is saturated or the model is not adequate. In a similar manner the implicit threshold value of concentrations may follow from mass conservation: the total mass of reactants, say $R$, must be constant. Therefore it makes sense to look for solutions $u(x,t)$ in the rectangle $\{u\in\R^N\mid 0\leq u_i\leq R_i,\; i=1,...,N\}$ or the simplex $\{u\in\R^N\mid\sum_{i=1}^Nu_i=R,\;u_i\geq 0,\;i=1,...,N\}$. In general equations of the form \eqref{problem 1} should be therefore considered under the presence of state constraints: $u(x,t)\in C$ for $x\in\Omega$, $t\in [0,T]$, where $C$ is a given closed subset of the phase space $\R^N$.\\
\indent In what follows we admit also discontinuous nonlinearities $f$ (or $g$).  This appears when, for instance, system data are determined by measurements or a subject to phase transition phenomena and is motivated by numerous applications of systems with hysteresis (see e.g. \cite{Terman}, \cite{Brokate}). The typical situation concerns \eqref{problem 1} with $N=M=1$ and is of the form
$$u_t=u_{xx}+H(u),$$
where  $H$ is the  hysteresis operator - see \cite{Visi}, \cite{KrasnoHyst}. In the simplest case $H$ is driven by the Heaviside function and maybe described via the related Nemytskii operator: given a threshold value $\alpha>0$, an input function $u:[0,T]\to\R$ with $u(0)\leq \alpha$, then the output $H(u)(t)=0$ if $u(s)<\alpha$ for all $s\in [0,t]$ and $H(u)(t)=1$ if $u(s)=\alpha$ for some $s\in [0,t]$. For some other instances of the problem -- see \cite{Bernal}, \cite{Carl} and numerous examples in \cite{Bothe}. The common way to overcome this obstacle is to replace the discontinuous $f$, or $g$ in \eqref{advection}, by an appropriate set-valued regularization $F$ or $G$ (introduced e.g. by Fillipov or Krasovski -- see \cite[Sect. 2.7]{Fil} or \cite[p. 101]{AuCel}) and instead of \eqref{problem 1} consider a problem
\begin{gather}\label{problem 2}u_t\in Lu+F(t,x,u,\nabla u)\;\; \text{or}\\
\label{problem 2bis}u_t\in Lu-\Gamma u+G(t,x,u)\end{gather}
subject to  initial and boundary conditions,
where $F:[0,T]\times\Omega\times C\times\R^{MN}\multi\R^N$ (or $G:[0,T]\times\Omega\times C\multi\R^N$) is an upper semicontinuous set-valued map with compact convex values.

\subsection{Zeros of set-valued maps}\label{zeros}
In the present paper we shall deal with the existence of steady state (stationary solutions) of state constrained autonomous problems related to \eqref{problem 2} or \eqref{problem 2bis}. This leads to the second objective of the present paper.\\
\indent Since 1941, when Kakutani showed  that every upper
semicontinuous set-valued self map $\vp$ of the closed ball $D$ in $\R^n$ admitting closed convex values has a fixed point, a lot of attention has been paid to the different aspects of the fixed point theory for set-valued maps (see e.g. \cite{Gorn-mon}).  In one direction the development has led to substantial
weakening in the assumption that the values of the mapping are subsets of
its domain. The idea is well-illustrated by the classical (single-valued) mean value theorem of Bolzano.
\begin{thm}\label{Bolzano} If $f:C:=[a,b]\to\R$ is continuous,  $f(a)f(b)\leq 0$ (for instance $f(a)\geq 0$ and $f(b)\leq 0$), then there is $\bar x\in C$ such that $f(\bar x)=0$.\end{thm}
\indent This important observation has been generalized by Poincar\'e in 1883 in his famous conjecture proved by Miranda \cite{Miranda}.
\begin{thm}\label{Miranda}{\em (Poincar\'e-Miranda)}
Let $C=\Pi_{k=1}^n[a_k,b_k]$ be an $n$-dimensional cube and let
$F^-_k:=\{x\in C\mid x_k=a_k\}$, $F^+_k:=\{x\in C\mid x_k=b_k\}$, $k=1,2,...,n$, denote the $k$-th face of $C$. Let $f=(f_1,...,f_n):C\to\R^n$ be continuous and suppose that for all $k=1,...,n$
\be\label{tangency 1}f_k(x)\begin{cases}\geq 0 &\text{for every}\;x\in F_k^-\\
\leq 0 &\text{for every}\;\;x\in F_k^+.\end{cases}\ee
Then $f$ has a zero, i.e., there is $\bar x\in C$ such that $f(\bar x)=0$. \end{thm}
Quite a complicated history of this and other similar results is well-described by Mawhin \cite{Mawhin 1} and \cite{Browder 1} (see also \cite{Mawhin 2}, \cite{Mawhin 3}). In the spirit of the above we have (see \cite{Mawhin 3}, \cite{Sch})
\begin{thm}\label{Ghezzo} Let $C=\{x\in\ell^2\mid |x_k|\leq k^{-1}\}$ be the Hilbert cube, let $f:C\to\ell^2$ be continuous and such that for all $k\in\N$
\be\label{tangency 2}f_k(x_1,...,x_{k-1},-\frac{1}{k},x_{k+1},...)\geq 0,\;\; f_k(x_1,...,x_{k-1},\frac{1}{k},x_{k+1},...)\leq 0,\ee
then $f$ has a zero.
\end{thm}
\indent In order to understand the nature of assumptions of these results we need to recall the following definition. Let $E$ be a Banach space, $K\subset E$ be a closed set and $x\in K$. The contingent (or Bouligand) cone $T_K(x)$ is defined by
$$T_K(x)=\{v\in E\mid \liminf_{h\to
0^+} d_K(x+hv)/h=0\}$$
and the Clarke tangent cone is defined by
$$C_K(x)=\{v\in E\mid\lim_{h\to 0^+,\,y\to x,\,y\in K}d_K(y+hv)/h=0\},
$$
where $d_K(u):=\inf_{y\in K}\|y-u\|$.
$T_K(x)$ and $C_K(x)$ are closed cones; additionally $C_K(x)$ is convex. In general $C_K(x)\subset T_K(x)$ and if $K$ is convex, then
\be\label{rowne}T_K(x)=C_K(x)=\cl S_K(x)\ee where
$$S_K(x):=\bigcup_{h>0}h(K-x).$$
Observe that if $x$ belongs to the interior of $K$, then $T_K(x)=E$. For examples and a detailed discussion see \cite{Aubin}.
\begin{rem}\label{examples} (1)  If $C\subset \R^n$ is a cube (as in Theorem \ref{Miranda}) and $x\in C$ is a boundary point, then $v\in T_C(x)$ if and only if $v_k\geq 0$ if $x_k=a_k$ and $v_k\leq 0$ if $x_k=b_k$. This implies that assumption \eqref{tangency 1} holds true if and only if $f(x)\in T_C(x)$, $x\in C$. Similarly, in the context of Theorem \ref{Ghezzo}, assumption \eqref{tangency 2} is satisfied if and only if $f(x)\in T_C(x)$ for all $x\in C$.\\
\indent (2) If $K=D(0,r)$ is a closed ball in a Hilbert space, then $T_K(x)=\{v\in E\mid \la x,v\ra\leq 0\}$ for $x$ in the boundary of $K$.
\end{rem}
\indent In the apparently independent stream of research, the best known equilibrium result is the following pioneering result of Browder \cite{Browder 2} (with some modification due to Halpern and Bergman \cite{Halpern 1,Halpern 2}) being, in the opinion of Aubin and Cellina (see \cite[p. 213, Chapter 5.2]{AuCel} and the discussion therein), `one of the most powerful theorems of nonlinear analysis'.

\begin{thm}\label{Browder}Assume that $K\subset E$ is  compact convex and  $\vp:K\multi E$ is upper semicontinuous with closed
convex values. If $\vp$ satisfies the {\bf weak tangency} condition with
respect to $K$,  i.e.
\be\label{weak tangency}\forall\,x\in
K\;\;\;\; \vp(x)\cap T_K(x)\neq \emptyset,\ee then $\vp$ has an equilibrium:
there is $\bar x\in K$ such that $0\in\vp(\bar x)$. If $\vp$ satisfies the so-called the {\bf weak inwardness (or outwardness)} condition, i.e.
\be\label{weak inward}\forall\,x\in K\;\;\;\; \vp(x)\cap (x+T_K(x))\neq\emptyset\ee
{\em (}or $\vp(x)\cap (x-T_K(x))\neq\emptyset$ for $x\in K${\em )},
then $\vp$ has a fixed point. \end{thm}
In view of Remark \ref{examples} it is evident that Theorem \ref{Browder} provides a far reaching generalization of Theorems \ref{Miranda} and \ref{Ghezzo}.

\begin{rem}
It is easy to see that if $K$ is convex and $0\in K$, then $T_K(x)\in x+T_K(x)$ for all $x\in K$. Hence in Theorem \ref{Ghezzo}
(and in \ref{Miranda} if $a_k\leq 0\leq b_k$, $k=1,...,n$) $f$ has a fixed point since $f(x)\in x+T_C(x)$, $x\in C$ and, therefore $f(x)-x=0$ for some $x\in C$. Similarly in Theorem \ref{Browder} if $0\in K$, then \eqref{weak tangency} implies \eqref{weak inward} and $\vp$ has fixed points.
\end{rem}

Two drawbacks of this result has to be pointed out. In order to get a decent tool to study existence of equilibria one needs to get rid of compactness and convexity in Theorem \ref{Browder}. The best known result in this direction is due to Deimling -- see \cite[Th.11.5]{Deimling 1}, \cite{Deimling 3}.
\begin{thm}\label{Deimling} Let $K$ be a closed bounded convex subset of a Banach space $E$ and let  an upper semicontinuous map $\vp:K\multi E$ with compact convex values be condensing with respect to the Kuratowski or Hausdorff measure of noncompactness. If $\vp$ is weakly inward, then $\vp$ has a fixed point.\end{thm}
It is interesting to observe that in Deimling's theorem there is no way to replace inwardness by outwardness condition although  it was possible in Theorem \ref{Browder}.\\
\indent In order to discuss a nonconvex version of Theorem \ref{Browder} one needs to understand which property of a set is a suitable substitute for convexity and what should be a suitable counterpart of tangency. This problem was addressed in \cite{Ben-Krysz} and discussed in \cite{Krysz}.\\
\indent After \cite{Ben-Krysz} we say that a closed $K\subset E$ is an $\mathcal L$-{\bf retract} if there is $\eps>0$, a continuous $r:U=B(K,\eta)\to K$, where $B(K,\eta):=\{x\in E\mid d_K(x)<\eta\}$ (\footnote{Clearly $B(K,\eta)=K+\eta B$, where $B=B(0,1)$ is the open unit ball in $E$.}) and a constant $L\geq 1$ such that
$$r(x)=x\; \text{for}\; x\in K\;\;\text{and}\;\;\|r(x)-x\|\leq Ld_K(x)\; \text{for}\; x\in U.$$
Therefore $K$ is an $\mathcal L$-retract whenever $K$ is a neighborhood retract in $E$ with retraction $r$ such that distance of $x\in U$ from $r(x)\in K$ may be controlled by the distance $d_K(x)$. The class of $\mathcal L$-retracts is large. Closed convex sets (in this case one can define $r$ on $E$ with  $L=1+\eps$, where $\eps>0$ is arbitrary), compact sets being bi-Lipschitz homeomorphic with closed convex sets, the so-called proximate retracts, Lipschitz retracts and epi-Lipschitz sets (in the sense of Rockafellar \cite{Rocka}) are $\mathcal L$-retracts.
\begin{rem}\label{Euler} (1) If $X$ is a topological space of finite type, i.e., such that the (singular with rational coefficients) cohomology groups $H^k(X;\Q)$, $k\geq 0$, groups are  finitely generated and vanish above some dimension, then  the {\em Euler characteristic} $\chi (X):= \sum_{k=0}^\infty (-1)^k \dim H^k (X,\Q)$
of $X$ is well-defined. \\
\indent (2) If $X$ is a neighborhood retract in $E$ and $f:X\to X$
is compact, then $f$ is a {\em Lefschetz map}, i.e. the
homomorphism $H^*(f)$ is a {\em Leray endomorphism} of
$H^*(X,\Q)$ and the generalized Lefschetz number $\Lambda (f)$ of $f$ is well-defined -- see \cite[Def. V.(2.1), (3.1). Th.
(5.1)]{Dugundji-Granas}. If $f$ is homotopic
to the identity $I_X$ on $X$, then $H^*(f) =
H^*(I_X)$ is the identity $H^{*}(X)$. This implies that $I_X$ is a
Lefschetz map; hence $H^* (X)$ is of finite type and the Euler
characteristic $\chi(X)$ is well-defined. Moreover, in this case
$\Lambda(f)$ is  equal to the ordinary Lefschetz number
$\lambda(f)=\lambda(I_X)=\chi(X)$ (for details concerning
these notions see also e.g. \cite{Brown}). In particular if $\chi(X)\neq 0$, then $f$ has a fixed point.
\end{rem}
In view of the above if $K$ is a compact $\mathcal L$-retract, then its Euler characteristic $\chi(K)$ is well-defined. Note that if $K$ is additionally convex, then $\chi(K)=1$. After \cite{Ben-Krysz} (see also \cite{cwisz-krysz 1}) we have the following result.
\begin{thm}\label{Ben-Krysz} Let $K\subset E$ be a compact $\mathcal L$-retract with
$\chi(K)\neq 0$. If $\vp:K\multi E$ is upper semicontinuous with closed
convex values and {\bf weakly tangent to $K$ in the sense of Clarke}, i.e. \be\label{ctan}\forall\,x\in
K\;\;\;\; \vp(x)\cap C_K(x)\neq\emptyset,\ee then $\vp$ has an
equilibrium.\end{thm}
Note that in condition \eqref{ctan} the Bouligand cone has been replaced by the Clarke cone; there are examples showing that \eqref{weak tangency} is not sufficient (see \cite{Krysz}); however if $\vp=f$ is a single-valued map, then \eqref{weak tangency} implies \eqref{ctan}. It is also evident that the weak inwardness in the sense of Clarke cone implies the existence of fixed points.\\
\indent There is no direct generalization of the equilibrium problem from Theorem \ref{Ben-Krysz} in the noncompact setting, although there were some partial
answers have been discussed in \cite{cwisz-krysz 1} and \cite{cwisz-krysz 2}, since we have the following example showing a compact tangent map without zeros.
\begin{ex} Let $E=\ell^2$ be the classical Hilbert space. Let $D_1=\{x\in E\mid \|x\|\leq 1/\sqrt{2}\}$ and let $f:D_1\to D_1$ be the variant of the famous Kakutani map: $f(x)=(\sqrt{1/2-\|x\|^2},x_1,x_2,...)$ for $x=(x_i)\in D_1$. Then $f$ is continuous, it has neither zeros nor fixed points and $\|f(x)\|=1/\sqrt{2}$ for all $x\in D_1$. Let $D$ be the unit closed ball in $E$ and $r:D\to D_1$ the radial retraction. Define $g:D\to D$ by
\[
g(x):=
\begin{cases}
-f(x)&\quad \text{for }x\in D_1;\\
\le 2\|x\|^2-2\pr f(r(x)) + \le 1 - 2\|x\|^2\pr x& \quad \text{for }x\in D\setminus D_1.
\end{cases}
\]
One can see that $g$ is well-defined, continuous and $g(x)=-x$ whenever $x\in \partial D$; an easy argument yields $g(x)\neq 0$ for every $x\in D$. Finally, define $\kappa :E\to E$ by
\[
\kappa(x_1,x_2,\ldots,x_n,\ldots):= \le x_1,2^{-1}x_2,\ldots,n^{-1}x_n,\ldots \pr.
\]
Clearly $\kappa$ is an injective compact linear map.
Thus $G:=\kappa\circ g :D\to D$ is compact and $G(x) \neq 0$ for  $x\in D$.
If $x\in \partial D$, then $\langle G(x),x \rangle = - \sum_{n=1}^\infty \frac{1}{n}x_n^2 \leq 0$; by Remark \ref{examples} (2), $G$ is tangent to $D$.
\end{ex}
For examples, further  generalizations and a deeper discussion of issues surveyed above the reader can see \cite{Krysz}.

 The main aim of the present paper is to show a result in this direction with applications to constrained steady state problems  related to \eqref{problem 2} or \eqref{problem 2bis}.

\section{Existence results}\label{third}

\subsection{The setting and results}\label{existence}  In order to study the existence of steady states of autonomous problem \eqref{problem 2} we shall take an appropriate appropriate abstract setting and consider the following coincidence problem
\be\label{AE} 0\in Au+\Phi(u),\;\; u\in K\subset E,\ee
where:
 \begin{enumerate}
 \item[$(A_1)$] $(E,\|\cdot\|)$ is a Banach space, $K\subset E$ is a closed convex set;
 \item[$(A_2)$] $A:E\supset D(A)\to E$ is a densely defined linear operator such that, for some $\omega\in\R$, $(\omega,\infty)\subset\rho(A)$;
 \item[$(A_3)$] there is another Banach space $(E_0,\|\cdot\|_0)$ and a closed convex $K_0\subset E_0$ such that $D(A)\subset E_0\subset E$, $K_0\subset K$ and  the identities $D(A)\hookrightarrow E_0$, $j:E_0\hookrightarrow E$ are continuous (\footnote{On $D(A)$ the {\em graph norm} $\|\cdot\|_A$ is considered: $$\|u\|_A:=\|u\|+\|Au\|,\;\; u\in D(A).$$});
\item[$(A_4)$] $\Phi:K_0\multi E$ is $H$-upper semicontinuous, bounded and has convex weakly compact values;
    \item[$(A_5)$] for all $u\in K_0$, $\Phi(u)\cap T_K(j(u))\neq\emptyset$.
\end{enumerate}

\begin{rem}\label{ass comments}
(1) Let us recall the so-called Lions construction, see \cite{Lions}, \cite{Arendt}. Let $V$ be a reflexive Banach space which is dense in a (real)  Hilbert space $H$ and suppose that the identity $V\to H$ is continuous. Suppose a bilinear continuous form $a:V\times V\to\R$ is such that
\be\label{Lions}\forall\,v\in V\;\; a(v,v)+\omega\|v\|_H^2\geq \alpha\|v\|_V^2,\ee
where $\omega\in\R$ and $\alpha>0$. Let $D(A)$ be the set of all $u\in V$ such that the function $V\ni v\mapsto a(u,v)$ is continuous on $V$ with the $H$-norm, i.e. such that there is $\beta_u\geq 0$ with $|a(u,v)|\leq\beta_u\|v\|_H$ for all $v\in V$. For any $u\in D(A)$ there is a unique $Au\in H$ such that
$$a(u,v)=-\la Au,v\ra_H,\;\; v\in V.$$
This defines a linear $A:D(A)\to H$. According to e.g. \cite[Prop. 4.1]{Show}, $D(A)$ is dense in $H$, $A$ is closed and $(\omega,+\infty)\subset\rho(A)$. Moreover, $A$ is maximal $\omega$-dissipative, i.e. it is the generator of a $C_0$-semigroup $\{S(t)\}_{t\geq 0}$ of linear operators on $H$ with growth bound equal to $\omega$. Putting $E:=H$, $E_0:=V$ and taking a closed convex bounded $K\subset E$ we see that assumption $(A_1)$ -- $(A_3)$ are satisfied provided $K_0:=K\cap E_0$. The situation described in this example is very typical in various applications.\\
\indent (2) Assumption $(A_4)$ is motivated by applications. $H$-upper semicontinuity (where $H$ stands for `Hausdorff') in $(A_4)$ means that for each $x\in K_0$ and $\eps>0$ there is $\delta>0$ such that $\Phi(y)\in B(\Phi(x),\eps)$ if $y\in K_0$ and $\|y-x\|<\delta$. It is well know that (see e.g. \cite[prop. 2.3]{Bothe}) that $\Phi$, being $H$-upper semicontinuous with weakly compact values, is upper semicontinuous when $E$ is endowed with the weak topology. This, in turn, implies, that given sequences $(x_n)\subset K_0$ and $y_n\in\Phi(x_n)$, if $x_n\to x_0\in K_0$, there is a subsequence $y_{n_k}$ such that $y_{n_k}\rightharpoonup y_0\in\Phi(x_0)$ ($\rightharpoonup$ denotes the weak convergence). Obviously if $\Phi$ is upper semicontinuous with closed compact values, then $(A_4)$ is satisfied, too.\\
\indent (3) Let $h>0$ and $h\omega<1$, then $h^{-1}\in\rho(A)$ and {\em the resolvent}
$$J_h:=(I-hA)^{-1}:E\to E$$
is well defined and continuous. Observe that $J_h(u)\in D(A)\subset E_0$ for any $u\in E$. Moreover for any $u\in E$
$$\|J_h(u)\|_A=\|J_h(u)\|+\|AJ_h(u)\|\leq ((1+h^{-1})\|J_h\|+h^{-1})\|u\|,$$
This, together with $(A_3)$, shows that $J_h:E\to E_0$ is continuous, too. Observe that $J_h$  is compact for some $h>0$ with $h\omega<1$ if and only if the identity  $D(A)\to E$ is compact. If $D(A)\to E_0$ is compact, then $J_h$, as a map from $E$ to $E_0$ is compact for all $h>0$ with $h\omega<1$. It is worth to note that in the situation described in part (1), $\|J_h\|\leq (1-h\omega)^{-1}$.\end{rem}

Let us recall a version of Lemma 17 from \cite{bader-krysz}; for the reader's convenience we give an independent proof.

\begin{lem}\label{convex selection} For any $\eps>0$, there exists a locally Lipschitz map $f=f_\eps:K_0\to E$ being an $\eps$-graph-approximation {\em (\footnote{The graph of $f$ is contained in the $\eps$-neighborhood of the graph of $\Phi$, i.e. $f(x)\in \Phi(B_0(x,\eps)\cap K_0)+\eps B$, where $B$ is the unit open ball in $E$, for any $x\in K_0$; in particular $f$ is bounded. Here and below we write $B(x,r)$, $x\in E$, to denote a ball in $E$ and $B_0(x,r)$, $x\in E_0$, to denote a ball in $E_0$.})}
of $\Phi$ and, for all $u\in K_0$,
\be\label{tangency approx}f(u)\in T_K(j(u)).\ee\end{lem}

\noindent {\sc Proof:} Take $\eps>0$ and $u\in K_0$. By $(A_5)$ and \eqref{rowne} there is $v(u)\in E$ such that
$$v(u)\in B(\Phi(u),\eps/4)\cap S_K(j(u)).$$
Hence, there is $\alpha(u)>0$ such that
$$j(u)+\alpha(u)v(u)\in K.$$
By the $H$-upper semicontinuity choose a number $\gamma(u)$, $0<\gamma(u)<\eps/4$ such that $\Phi(B_0(u,2\gamma(u))\cap K_0)\subset B(\Phi(u),\eps/2)$ and a number
$0<\delta(u)<\min\{\gamma(u)/C,\gamma(u)/\alpha(u)\}$, where $C:=\|j\|$.\\
\indent Let $\{\lambda_s\}_{s\in S}$ be a locally finite locally Lipschitzian
partition of unity refining the open cover $\{B_0(u,\delta(u)\alpha(u))\cap K_0\}_{u\in
K_0}$. For any $s\in S$, there is $u_s\in K$ such that the support $\supp\lambda_s\subset
B_0(u_s,\delta_s\alpha_s)\cap K_0$ where we have put $\delta_s:=\delta(u_s)$ and
$\alpha_s:=\alpha(u_s)$. Additionally let us set
$v_s:=v(u_s)$ and $\gamma_s:=\gamma(u_s)$.\\
\indent For any $s\in S$, we define a map $f_s:K_0\to E$ by the formula
$$f_s(u):=\frac{1}{\alpha_s}(j(u_s)-j(u))+v_s,\;\;\;u\in K_0.$$
Observe, that for $s\in S$, $u\in K_0$, $j(u)+\alpha_sf_s(u)=j(u_s)+\alpha_sv_s\in K.$
Hence, for all $u\in K_0$,
$$f_s(u)\in S_K(j(u))\subset T_K(j(u)).$$
It is clear that $f_s$, $s\in S$, is Lipschitz continuous.\\
\indent Now we define $f:K_0\to E$ by the formula
$$f(u):=\sum_{s\in S}\lambda_s(u)f_s(u),\;\; u\in K_0.$$
Observe that $f$ is locally Lipschitz because so are all functions $\lambda_s$, $f_s$ for $s\in S$, and the covering $\{\supp\lambda_s\}_{s\in S}$ is
locally finite. Moreover, since, for $u\in K_0$, $f(u)$ is a (finite) convex
combination of vectors $f_s(u)\in T_K(j(u))$ and since $T_K(j(u))$ is convex, we see that $f(u)\in T_K(j(u))$ for all $u\in K_0$.\\
\indent Take  $u\in K_0$ and let $S(u)=\{s\in S\mid u\in \supp\lambda_s\}$. It
is clear that $S(u)$ is a finite set and
$$f(u)=\sum_{s\in S(u)}\lambda_s(u)f_s(u).$$
For any $s\in S(u)$, we have $u\in\supp\lambda_s\subset
B_0(u_s,\delta_s\alpha_s)\cap K_0$, i.e.
$$\|u-u_s\|_0<\delta_s\alpha_s<\gamma_s\;\;\;\mbox{
and }\;\;\; \|f_s(u)-v_s\|<\delta_s C<\gamma_s.$$ There is $s_0\in S(u)$ such
that $\gamma_{s_0}=\max_{s\in S(u)}\gamma_s$. If $s\in S(u)$, then
$$\|u_s-u_{s_0}\|_0\leq\|u_s-u\|_0+\|u_{s_0}-u\|_0<\gamma_s+\gamma_{s_0}\leq 2\gamma_{s_0}.$$
Therefore, for any $s\in S(u)$, \begin{gather*} f_s(u)\in
B(v_s,\gamma_{s_0})\subset \Phi(u_s) +(\eps/4+\gamma_{s_0})B\subset
\Phi(B_0(u_{s_0},2\gamma_{s_0})\cap K_0) +(\eps/4+\gamma_{s_0})B\\\subset
B(\Phi(u_{s_0}),\eps/4+\eps/2+\gamma_{s_0})\subset
\Phi(u_{s_0})+\eps B.\end{gather*} Hence, by convexity of
$\Phi(u_{s_0})+\eps B$,
$$f(u)\in \Phi(u_{s_0})+\eps B\subset
\Phi(B_0(u,\gamma_{s_0})\cap K_0)+\eps B\subset \Phi(B_0(u,\eps)\cap K_0)+\eps B.\eqno\square$$

\begin{lem}\label{lem.stycznosc} For every $u\in K_0$, we have
\be\label{wzor.stycznosc}
\lim_{h\to 0^+,\,v\to u,\, v\in K_0}\frac{1}{h}d(j(v)+hf(v),K)=0,
\ee where $f$ comes from Lemma \ref{convex selection}.\end{lem}

\noindent {\sc Proof:} Choose $u\in K_0$ and $\eps>0$.
Taking into account \eqref{tangency approx}, \eqref{rowne} and the continuity of $f$ there is $\delta>0$ such that if $v\in K_0$, $\|v-u\|_0<\delta$ and $0<h<\delta$ then
\[
d(j(v)+hf(u);K)<\eps h/2\quad\text{and}\quad \|f(v)-f(u)\| <\eps/2.
\]
Thus, for such $v$ and $h$ we have
$$
d(j(v)+hf(v),K) \leq d(j(v)+hf(u),K) + h\|f(u)-f(v)\|< \eps h.\eqno\square
$$

Now we are ready to prove

\begin{thm}\label{main 1} In addition to $(A_1)$ -- $(A_5)$ above, let us assume that $K$ is bounded and for all sufficiently small $h>0$:
\begin{enumerate}
\item[$(A_6)$] $J_h(K)\subset K_0$;
\item[$(A_7)$] $J_h:E\to E_0$ is compact.
\end{enumerate}
Then there is $u\in K_0\cap D(A)$ such that $0\in Au+\Phi(u)$.\end{thm}

\noindent {\sc Proof:}  Choose $\eps>0$ and $f=f_\eps$ according to Lemma \ref{convex selection}. Denote by $r: E\to K$ an $\mathcal L$-retraction onto  $K$, i.e. $\|r(u)-u\| \leq Ld_K(u)$ for $u\in E$ for some $L\geq 1$.
For every $h>0$ with $h\omega<1$, the map $J_h\circ r(j + hf):K_0 \to K_0$ is well-defined due to $(A_6)$. Moreover $J_h\circ r(j+ hf)$ is continuous and compact, since $K$ and $\Phi$ (and so does $f$) are bounded and $J_h$ is compact. Then, by the Schauder fixed point theorem, for large $n\geq 1$ (precisely for $n>\omega$), there is $u_n$ in $K_0$ such that
\[
u_n = J_{1/n}\circ r\left(j(u_n) +\frac{1}{n}f(u_n)\right),
\]
so $u_n\in D(A)$ and
\[
j(u_n)-\frac{1}{n}Au_n = r\left(j(u_n) +\frac{1}{n}f(u_n)\right).
\]
As a result, we have
\[
Au_n + f(u_n) =  n\left(j(u_n) + \frac{1}{n}f(u_n) - r\left(j(u_n) +\frac{1}{n}f(u_n)\right)\right).
\]
Therefore
\be\label{dowod.szacowanie.stycz}
\|Au_n + f(u_n)\|\leq nL d_K\left(j(u_n)+\frac{1}{n}f(u_n)\right)\leq L\|f(u_n)\|.\ee
Hence $\{ Au_n\}_{n\geq 1}$ is bounded in $E$. Fix $h>0$ with $h\omega<1$ and note that
\[
\{u_n\}_{n\geq 1} = J_h\left(\left\{ j(u_n) - hAu_n \right \}_{n\geq 1}\right).
\]
Since $\{ j(u_n) - hAu_n\}_{n\geq 1}$ is bounded in $E$, the above equality yields $\left \{ u_n \right \}_{n\geq 1}$ is relatively compact in $E_0$.
Passing to a subsequence if necessary, we can assume that $u_n\to u_\eps$ in $E_0$ and $u_\eps\in K_0$.
In view of \eqref{dowod.szacowanie.stycz} and Lemma \ref{lem.stycznosc} we have
\be\label{koniec}
\|Au_n + f(u_n)\|= n d_K\left(j(u_n)+\frac{1}{n}f(u_n)\right) \to 0,\; \text{ as }n\to\infty.
\ee
Hence $Au_n \to -f(u_\eps)$ in $\E$.
The closedness of $A$ yields $u_\eps\in D(A)\cap K_0$ and $-Au_\eps=f(u_\eps)$.\\
\indent Arguing as above, we may assume without loss of generality  that $u_\eps\to u_0\in K_0$ as $\eps\to 0$. Let $v_\eps:=-Au_\eps$. Since
$$v_\eps:=f(u_\eps)\in \Phi(B_0(u_\eps,\eps)\cap K_0)+\eps B,$$
there is $u_\eps'\in K_0$ and $v_\eps'\in\Phi(u_\eps')$ such that $\|u_\eps-u_\eps'\|_0<\eps$ and $\|v_\eps-v_\eps'\|<\eps$. Clearly $(u_\eps',v_\eps')\in \Gr(\Phi)$ and $u_\eps'\to u_0$; in view of Remark \ref{ass comments} (2) we gather that, after passing to a subsequence if necessary, $v_\eps'\rightharpoonup v_0\in\Phi(u_0)$. Since $v_\eps\rightharpoonup v_0$, too, and the graph of $A$, being closed and convex, is also weakly closed, we see that $u_0\in D(A)$ and $-Au_0=v_0\in\Phi(u_0)$.\hfill $\square$

Now we are going to establish a counterpart of Theorem \ref{main 1} valid for $\mathcal L$-retracts. In this case the choice of $E_0$ is immaterial since we shall assume that $\Phi$ is defined on $K$. In addition Let us assume that:
\begin{enumerate}
\item[$(B_1)$] $E$ is a Banach space and $K\subset E$ is a bounded $\mathcal L$-retract;
\item[$(B_2)$] $A:D(A)\to E$ is a densely defined linear operator such that for some $\omega\in\R$, $(\omega,+\infty)\subset\rho(A)$ and $\|(A-\lambda I)^{-1}\|\leq (\lambda-\omega)^{-1}$ for $\lambda>\max\{0,\omega\}$ (\footnote{I.e. $A$ is maximal $\omega$-dissipative.});
\item[$(B_3)$] $\Phi:K\multi E$ is bounded $H$-upper semicontinuous with convex weakly compact values;
\item[$(B_4)$] $K$ is resolvent invariant, i.e. for $h>0$ with $h\omega<1$, $J_h:E\to E$ is compact and $J_h(K)\subset K$ for sufficiently small $h$.
    \end{enumerate}

\begin{rem}\label{Euler 2} (1) In view $(B_2)$ and $(B_4)$, the Euler characteristic $\chi(K)$ is well-defined. Indeed, by Remark \ref{Euler}, it is sufficient to prove that the identity $I_K$ is homotopic to a compact map. To this end fix  $h>0$ with $h\omega<1$, such that $J_h$ is compact and  consider $h:K\times [0,1]\to K$ given by the formula
$$h(x,t)=\begin{cases}J_{th}(x)&\text{if}\;\;\; t\in (0,1],\\
x&\text{if}\;\;\; t=0,\end{cases}\;\;x\in K.$$
Assumption $(B_2)$ implies that $\lim_{t\to 0^+}J_{th}(x)=x$ for all $x\in E$. Moreover the map $E\times (0,1]\ni (x,t)\mapsto J_{th}(x)\in E$ is continuous. Thus $h$ is continuous and provides a homotopy joining the identity to the compact map $J_h$. As a consequence if $\chi(K)\neq 0$, then any compact map $g:K\to K$ homotopic to the identity has fixed points.\\
\indent (2) If $A$ is given as in Remark \ref{ass comments} (with $E=H$), then assumption $(B_2)$ is satisfied. Moreover, in this case $(B_4)$ holds for a convex $K$ if and only if $K$ is semigroup invariant, i.e. $S(t)K\subset K$ for any $t\geq 0$. Indeed if $K$ is resolvent invariant then by the Post-Widder formula (see \cite[Corollary 5.5, 5.6]{Engel}) for each $x\in K$ and $t>0$
$$ S(t)x=\lim_{n\to\infty}J_{t/n}^nx\in K.$$
Conversely, if $K$ is semigroup invariant, then by \cite[Th. 1.10]{Engel}, for any $h>0$ with $h\omega<1$ and $x\in K$,
$$J_hx=\frac{1}{h}\int_0^\infty e^{-t/h}S(t)x\,dt\in K.$$
\end{rem}

\begin{thm}\label{main 2} Under assumptions $(B_1)$ -- $(B_4)$, the problem \eqref{AE} has a solution if the Euler characteristic $\chi(K)\neq 0$ and $\Phi$ satisfies the weak tangency condition in the sense of Clarke cones, i.e.
\be\label{Clarke tangency}\forall\,u\in K\;\;\;\; \Phi(u)\cap C_K(u)\neq\emptyset.\ee
\end{thm}
First we need a result, which may be of interest on its own, similar to Lemma \ref{convex selection}.
\begin{lem}\label{selection} Suppose $X\subset E$ is closed and that $\Phi:X\multi E$ is $H$-upper semicontinuous with convex values. Let a function $\xi:X\times E\to\R$ be such that for each $z\in E$, $\xi(\cdot,z)$ is upper semicontinuous (as a real function) and for each $x\in X$, $\xi(x,\cdot)$ is convex. If for all $x\in X$, $\inf_{z\in\Phi(x)}\xi(x,z)\leq 0$, then for any $\eps>0$ there exists a locally Lipschitz $\eps$-graph-approximation $f=f_\eps:X\to E$ of $\Phi$ such
$\xi(x,f(x))<\eps$ for all $x\in X$.\end{lem}

\noindent {\sc Proof:} For any $z\in X$ choose $0<\delta_z<\eps$ such that $\Phi(B(z,\delta_z)\cap X)\subset \Phi(z)+\eps B$ and let an open covering $\mathcal U$ of $X$ be a star refinement of the covering $\{B(z,\delta_z)\cap X\}_{z\in X}$ of $X$. \\
\indent For each $x\in X$ choose $z_x\in\Phi(x)$ such that $\xi(x,z_x)<\eps$. Given $U\in {\mathcal U}$ and $x\in U$ let
$$V_U(x):=\{y\in U\mid \xi(y,z_x)<\eps\}.$$
Clearly $x\in V_U(x)$. Hence ${\mathcal V}:=\{V_U(x)\}_{U\in {\mathcal U},\,x\in U}$ is an open cover of $X$. Let $\{\lambda_s\}_{s\in S}$ be a locally Lipschitz partition of unity subordinated to $\mathcal V$, i.e. for any $s\in S$, there is $U_s\in {\mathcal U}$, $x_s\in U_s$ such that $\supp\lambda_s\subset V_s:=V_{U_s}(x_s)$. Let
$$f(x):=\sum_{s\in S}\lambda_s(x)z_s,\;\; x\in X,$$
where $z_s:=z_{x_s}$. Then $f$ is well-defined and locally Lipschitz.\\
\indent Take $x\in X$ and let $S(x):=\{s\in S\mid\lambda_s(x)\neq 0\}$. If $s\in S(x)$, then $x\in V_s$, i.e. $\xi(x,z_s)<\eps$. By convexity of $\xi(x,\cdot)$ we gather that $\xi(x,f(x))<\eps$. Since $x_s\in U_s$ and $\mathcal U$ is a star refinement of $\{B(z,\delta_z)\cap X\}$ we get that for all $s\in S(x)$,
$x, x_s$ belong to the star of $x$ with respect to $\mathcal U$:
$$x,x_s\in\bigcup_{\{U\in {\mathcal U}\mid x\in U\}}U\subset B(z,\delta_z)\cap X$$
for some $z\in Z$.  Hence $z\in B(x,\eps)$ and  for $s\in S(x)$, $z_s\in\Phi(x_s)\subset \Phi(z)+\eps B$. This together with the convexity of $\Phi(z)$ shows that
$$f(x)\in \Phi(z)+\eps B\subset \Phi(B(x,\eps))+\eps B.\eqno\square$$

\noindent {\sc Proof} of Theorem \ref{main 2}: Let
$$\xi(u,v):=\sup_{p\in\partial d_K(u)}\la p,v\ra \in\R,\;\; u\in K,\;v\in E,$$
where $\partial d_K(u)\subset E^*$ denotes the generalized Clarke gradient at $u\in K$ of the (locally Lipschitz) function $d_K$. It is clear that $$\xi(u,v)=d_K^\circ(u;v)=\limsup_{y\to u,\,y\in K,\,t\to 0^+}\frac{d_K(y+tv)}{t}$$
is the Clarke directional derivative of $d_K$ at $u$ in the direction $v$. Then $\xi:K\times E\to\R$ is upper semicontinuous and, for each $u\in K$, $\xi(u,\cdot)$ is convex.\\
\indent Observe now that
$$C_K(u)=[\partial d_K(u)]^-:=\{v\in E\mid \xi(u,v)\leq 0\}.$$
Condition \eqref{Clarke tangency} together with with $(B_3)$ show that all assumptions of Lemma \ref{selection} are satisfied.\\
\indent Take an $\mathcal L$-retraction $r:B(K,\eta)\to K$ with constant $L$. Since $\Phi$ is bounded, there are $\lambda_0>0$ and $\eps_0>0$ such that
$$\forall\,u\in K\;\;\; u+\lambda f(u)\in B(K,\eta),$$
for any $f:K\to E$ being an $\eps$-graph-approximation of $\Phi$ with $0<\eps<\eps_0$ and $0<\lambda\leq\lambda_0$.\\
\indent Suppose now to the contrary that there are no solutions to \eqref{AE}. We claim that there is $0<\eps<\eps_0$ such that if $u\in K\cap D(A)$ and $f$ is an $\eps$-graph-approximation of $\Phi$, then
\be\label{war 1}\|Au+f(u)\|\geq (L+1)\eps.\ee
If not then there are sequences $\eps_0>\eps_n\to 0^+$, $u_n\in K$ and an $\eps_n$-approximation $f_n:K\to E$ of $\Phi$ such that
$$\|Au_n-f_n(u_n)\|<(L+1)\eps_n,\;\; n\in\N.$$
This implies that the sequence $(Au_n)$ is bounded; hence by the same argument as in the proof of Theorem \ref{main 1}, we gather that, passing to a subsequence if necessary, $u_n\to u_0\in K$.\\
\indent Since $f_n(u_n)\in\Phi(B(u_n,\eps_n))+\eps_nB$, we find $u_n'\in B(u_n,\eps_n)$ and $v_n'\in\Phi(u_n')$ such that $\|f_n(u_n)-v_n'\|<\eps_n$. By Remark \ref{ass comments} (2) we may assume that $v_n'\rightharpoonup v_0\in\Phi(u_0)$. This implies that $f_n(u_n)\rightharpoonup v_0$ and, thus $-Au_n\rightharpoonup v_0$, too. Hence $v_0= -Au_0$, i.e. $0\in Au_0+\Phi(u_0)$: a contradiction.\\
\indent Now take $\eps>0$ provided above and, using Lemma \ref{selection}, let $f:K\to\E$ be an $\eps$-graph-approximation of $\Phi$ such that $\xi(u,f(u))<\eps$ for all $u\in K$. Take a decreasing sequence $h_n\to 0^+$ with $h_1<\lambda$. Since $f$ is bounded the map
$$K\ni u\mapsto g_n(u):=J_{h_n}\circ r(u +h_nf(u))\in K$$
is well-defined a compact. Moreover $h:K\times [0,1]\to K$ given by
$$h(u,t):=\begin{cases}J_{th_n}\circ r\left(u +th_nf(u)\right)&\text{if}\;\;\; t\in (0,1];\\
u&\text{if}\;\;\; t=0,\end{cases}\;\;u\in K$$
provides a (continuous) homotopy joining the identity on $K$ with $g_n$. In view of Remark \ref{Euler 2}, $g_n(u_n)=u_n$ for some $u_n\in K\cap D(A)$. This means that
\be\label{war 2}Au_n+f(u_n)= h_n^{-1}(u_n+h_nf(u_n)-r(u_n+h_nf(u_n))).\ee
 Similarly as before we may suppose that $u_n\to u_0\in K$; therefore $f(u_n)\to f(u_0)$. By \eqref{war 1} and \eqref{war 2}, for all $n\in\N$,
$$(L+1)\eps\leq \|Au_n+f(u_n)\|\leq Lh_n^{-1}d_K(u_n+h_nf(u_n))\leq Lh_n^{-1}d_K(u_n+h_nf(u_0))+L\|f(u_n)-f(u_0)\|.$$
Passing to $\limsup$ and remembering that $\xi(u_0,f(u_0))<\eps$ we get
$$(L+1)\eps\leq\limsup_{n\to\infty}\|Au_n+f(u_n)\|\leq L\eps,$$
a contradiction. This completes the proof.\hfill $\square$

\subsection{Invariance and viability}\label{invariance} A central role among assumptions of Theorems \ref{main 1} and \ref{main 2} is played by the resolvent invariance of the set $K$ and the tangency condition. Let us consider conditions $(B_1)$ -- $(B_4)$ and \eqref{Clarke tangency} and let $K$ be an arbitrary closed subset of $E$.  The Hille-Yosida Theorem implies that in this case $A$ is the generator of a $C_0$ semigroup $\{S(t)\}_{t\geq 0}$. It is not difficult to show that $(B_4)$ and \eqref{Clarke tangency} imply that
$$\forall\,u\in K\;\;\; \liminf_{v\to u,\,v\in K,\,h\to 0^+}\frac{\dist(J_h(v+h\Phi(v)),K)}{h}=0\;\; (\footnote{Here $\dist$ stands for the distance between sets, i.e. $\dist(X,Y):=\inf\{d(x,y)\mid x\in X,\;y\in Y\}$.}).$$
This condition implies that
\be\label{strong}\forall\,u\in K\;\;\; \Phi(u)\cap T_K^A(u)\neq\emptyset,\ee
where
$$T_K^A(u)=\left\{v\in E\mid \liminf_{h\to 0^+}\frac{d_K(S_v(h)u)}{h}=0 \right\}$$
and $[0,+\infty)\ni t\mapsto S_v(t)u$ is the mild solution to the problem $\dot x=Ax+v$, $x(0)=u$. Finally \eqref{strong} is equivalent to the following:
\begin{eqnarray}\label{npm}\hbox{the problem}\;\; \dot
u\in Au+\Phi(u),\; u(0)=x\in M\;\;\hbox{has}\\\nonumber
\hbox{a mild solution}\;\; u:[0,+\infty)\to E\;\;
\hbox{takie, że}\;\;u(t)\in M\;\;\hbox{dla}\;\; t\geq
0,\end{eqnarray} because the semigroup $\{S(t)\}_{t\geq 0}$ is immediately compact. This and related results are thoroughly discussed in \cite{Vrabie 1} and \cite{Bothe}. We thus see that our conditions imply the {\bf invariance} of $K$ (sometimes called viability) with respect to the `heat flow' generated by $A$, i.e. condition \eqref{npm}. Conversely condition
\be\label{strong 1}\forall\,u\in K\;\;\; 0\in T_K^A(u)\ee
implies the semigroup invariance and, in case of a convex $K$, resolvent invariance $(B_4)$. The point is that, in concrete situations of differential problems, condition \eqref{strong 1} needs to be verified. In most cases this can be done via an appropriate use of the maximum principles. In the next section we shall encounter examples of such arguments.\\
\indent The problem of invariance of systems of parabolic PDE was studied in numerous papers \cite{Amann}, \cite{Kuiper}, \cite{Smoller}, \cite{bader-krysz}, \cite{Bothe}, \cite{Wein} (and references therein), \cite{Walter} (the so-called M\"uller conditions important in various applications). The most general, often necessary and sufficient, abstract results are presented in \cite{Vrabie 1}. The invariance problem of parabolic problem from \eqref{npm} will  be studied in  the forthcoming paper \cite{krysz-siem}. In particular we shall study the topological structure of the set of all viable (i.e. `surviving' in $K$) solutions and show  its  relation with the existence of steady states, i.e. solutions to \eqref{AE}.

\section{Applications}\label{fourth}

\subsection{The Neumann problem I} We now study the existence of steady state solutions to \eqref{problem 2} and consider the problem
\be\label{problem II}\begin{cases}-Lu\in F(x,u,\nabla u),\;\;u(x)\in C\;\; \text{a.e. on}\;\; \Omega,\\
\mbox{}\;\; \frac{\partial}{\partial {\bf n}}u=0\;\; \text{on}\;\; \partial\Omega,\end{cases}\ee
where $C\subset\R^N$ is a compact and convex set, $\frac{\partial}{\partial {\bf n}}u=\left(\frac{\partial u_1}{\partial {\bf n}},....\frac{\partial u_N}{\partial {\bf n}}\right)$ denotes the outward normal derivative of $u$. We are going to find a {\em strong solutions}: a function $u\in H^2(\Omega,\R^N)$ such that $-Lu(x)\in F(x,u(x),\nabla u(x))$ for a.a. $x\in\Omega$ and $\left.\frac{\partial u_i}{\partial{\bf n}}\right|_{\partial\Omega}=0$, $i=1,...,N$, in the sense of trace.\\
\indent Let us make the following assumptions:
\begin{enumerate}
\item[$(D)$] For all $i=1,...,N$, $d_i=d\in C^1(\bar\Omega)$ and $d>0$;
     \item[$(F_1)$] $F:\Omega\times C\times \R^{MN}\multi\R^N$ is upper semicontinuous with compact convex values;
    \item[$(F_2)$] there is a nonegative $b\in L^2(\Omega)$ such that
    $\sup_{y\in F(x,u,v)}|y|\leq b(x)$ for a.a. $x\in \Omega$, $u\in C$ and $v\in\R^{MN}$;
    \item[$(F_3)$] $F$ is weakly tangent to $C$ with respect to the second variable, i.e. $F(x,u,v)\cap T_C(u)\neq\emptyset$ for all $x\in \Omega$, $u\in C$ and all $v\in\R^{MN}$.
        \end{enumerate}
We now put
$$\begin{gathered}E:=L^2(\Omega,\R^N),\;K:=\{u\in E\mid u(x)\in C\;\;\text{a.e. on}\;\;\Omega\};\\
E_0:=H^1(\Omega,\R^N),\;\; K_0:=\{u\in E_0\mid u(x)\in C\;\;\text{a.e. on}\;\;\Omega\};\\
D(A):=\left\{u\in H^2(\Omega,\R^N)\mid \left.\frac{\partial u}{\partial{\bf n}}\right|_{\partial\Omega}=0\right\},\;\;Au:=Lu,\;\,u\in D(A).\end{gathered} $$
Clearly assumptions $(A_1)$, $(A_2)$ (with $\omega=0$) and $(A_3)$ are satisfied, $K$ is closed convex and bounded; condition $(A_7)$ holds true since the embedding $D(A)\to E_0$ is compact.\\
\indent For any $u\in K_0$, let
\be\label{def phi}\Phi(u):=\{v\in E\mid v(x)\in F(x,u(x),\nabla u(x))\;\;\text{for a.a.}\;\;x\in\Omega\}.\ee
Evidently values $\Phi:K_0\multi E$ are nonempty and convex.

\begin{prop}\label{properties of Phi} The map $\Phi$ satisfies conditions $(A_4)$ and $(A_5)$.
\end{prop}
\noindent {\sc Proof:} It is straightforward to show that $\Phi(u)$ is weakly compact (we work in a Hilbert space, thus closed convex and bounded sets are convex weakly compact). Below we shall  prove a slight generalization of Proposition 6.2 from \cite{Bothe}. It implies immediately that $\Phi$ is $H$-upper semicontinuous.
\begin{lem}\label{H-usc} If $\psi:\Omega\times\R^d\to\R^N$ is upper semicontinuous with convex compact values and $\sup_{y\in\psi(x,u)}|y|\leq b(x)+a|u|$ for all $u\in\R^d$ and a.a. $x\in\Omega$, where $b\in L^2(\Omega)$ and $a>0$, then the Nemytskii operator
$$\Psi:{\mathcal E}:=L^2(\Omega,\R^d)\multi E,\; \Psi(u):=\{v\in E\mid v(x)\in\psi(x,u(x))\;\;\text{for a.a.}\;\; x\in\Omega\},\;\; u\in {\mathcal E},$$
is $H$-upper semicontinuous\end{lem}
\noindent {\sc Proof:} Suppose it is not the case: there are $\e _0>0$, a sequences $u_n \to  u_0$ in ${\mathcal E}$ and $v_n \in \Psi\le u_n \pr $ such that
\begin{equation}\label{w-k_do_sprzecznosci}
v_n \notin \Psi \le u_0 \pr + B_E(0,\e_0),\;\; n\geq 1.
\end{equation}
Up to a subsequence $\le u_n \pr _{n\geq 1}$ converges a.e. on $\Omega$ to $u_0$ and there is $h\in L^2 \le \Omega,\R \pr $ such that $|u_n \le x \pr | \leq h \le x \pr $ for a.e. $x\in \Omega$ and every $n\geq 0$. By assumption
\[
|v_n \le x \pr| \leq b \le x \pr + a|u_n \le x \pr |  \leq b \le x \pr + a h \le x \pr \text{ for }n\geq 0\text{ and }a.e.\; x\in \Omega.
\]
There is $\eta >0$ such that for $A \subset \Omega$ with Lebesgue measure $\mu \le A \pr <\eta$
\begin{equation}\label{wzor-absolutna-ciaglosc-calki}
\int _A  4\le b\le x \pr + a h\le x \pr \pr ^2\,dx  < \e_0^2 / 2.
\end{equation}
For each $n\geq 0$, the set-valued map $H_n:=\psi\le\cdot, u_n \le \cdot \pr \pr :\Omega \m \R^N$ is measurable and
if $w:\Omega \to \R^N$ is a measurable selection of $H_n$, then $w \in E$ since \begin{equation}\label{w-k_wzrostu}
|w\le x \pr | \leq b\le x \pr + a h \le x \pr  \quad \text{for a.e. }x\in \Omega.
\end{equation}
\indent By the Egorov and Lusin theorems (see \cite[Th. 1]{Averna} for a multivalued version of the Lusin theorem) there is a compact  $\Omega_\eta \subset \Omega$ such that $\mu \le \Omega \setminus \Omega_\eta \pr < \eta$, $u_n \to u_0$ uniformly on $\Omega_\eta$,
the restriction $u_0|_{\Omega_\eta}:\Omega_\eta \to \R^N$ is continuous and  $H_0|_{\Omega_\eta} : \Omega_\eta \m \R^N$ is $H$-lower semicontinuous.\\
\indent Let $\delta := \e_0 /\sqrt{2\mu \le \Omega \pr }$. We will show that there is $n_0$ such that if $n\geq n_0 $ and $x\in \Omega_\eta$, then
\[
H_n \le x \pr \subset H_0 \le x \pr + B_{\R^N} \le 0, \delta \pr .
\]
Suppose to the contrary that  there is a subsequence $\le n_j \pr _{j\geq 1}$ and a sequence $\le x_j \pr _{j\geq 1 }$ in $\Omega_\eta$ such that
\begin{equation}\label{w-k_H_0}
H_{n_j} \le x_j \pr \not\subset H_0 \le x_j \pr + B_{\R^N} \le 0 ,\delta \pr .
\end{equation}
We can assume that $x_j \to x_0 \in \Omega_\eta$, since $\Omega_\eta$ is compact. The continuity of $u_0|_{\Omega _\eta}$ and the uniform convergence $u_n \to u_0$ on $\Omega_\eta$ imply that $u_{n_j}(x_j)\to u_0(x_0)$ and thus $\le x_j, u_{n_j}\le x_j \pr \pr \to \le x_0, u_0 \le x_0 \pr\pr$ as $j\to \infty$. The upper semicontinuity of $\f$ together with the $H$-lower semicontinuity of $H_0$ on $\Omega _\eta $ show that
$H_{n_j}\le x_j \pr \subset H_0\le x_j \pr +B_{\R^N} \le 0 ,\delta\pr$
for sufficiently large $j$, which contradicts \eqref{w-k_H_0}.\\
\indent Let us fix $n\geq n_0$. For a.e. $x\in \Omega_\eta$ we have
\begin{equation}\label{w-k_selekcja_v_n}
v_n\le x\pr \in H_n \le x \pr \subset H_0\le x \pr + B_{\R^N} \le 0, \delta \pr.
\end{equation}
Observe that the map $\Omega_\eta \ni x\mm B_{\R^N} \le v_n \le x \pr ,\delta \pr \cap H_0\le x\pr$ is measurable and has nonempty values for a.e. $x\in \Omega_\eta$.
By the Kuratowski--Ryll-Nardzewski theorem, there is a measurable selection $v: \Omega_\eta \to \R^N $, i.e. $v\le x \pr \in B_{\R^N} \le v_n \le x \pr ,\delta \pr \cap H_0\le x\pr $ for a.e. $x\in \Omega_\eta$.
Clearly $v\in L^2 \le \Omega_\eta , \R^N\pr$ and for a.e. $x\in \Omega_\eta$, $|v_n\le x\pr -v\le x \pr| <\delta$. Thus
\[
\int _{\Omega_\eta} |v_n\le x \pr - v\le x \pr|^2\,dx < \delta^2 \mu \le \Omega_\eta \pr < \e_0^2 / 2.
\]
Take an arbitrary selection $w:\Omega \to \R^N$ of $H_0$, i.e. $w\le x\pr \in H_0 \le x\pr $ for a.e. $x\in \Omega$.
Let $\chi=\chi_{\Omega_\eta}$ be the indicator of $\Omega_\eta$.
Notice that $\chi v +(1-\chi)w:\Omega \to \R^N$ is a square-integrable selection of $H_0$ (we identify $v:\Omega_\eta \to \R^N$ with the function $v:\Omega \to \R^N$ putting $v\equiv 0$ on $\Omega\setminus \Omega_\eta$).
By \eqref{w-k_wzrostu}
\[
|v_n\le x\pr - w \le x \pr | \leq |v_n\le x\pr | + |w\le x\pr |\leq 2\le b\le x \pr  + a h \le x \pr \pr\quad\text{for a.e. }x\in \Omega\setminus \Omega_\eta.
\]
Recall that  $\mu \le \Omega \setminus \Omega _\eta\pr<\eta$, hence and by \eqref{wzor-absolutna-ciaglosc-calki}
\begin{multline*}
\| v_n - \chi v+ (1- \chi)w\|^2= \int _{\Omega_\eta} |v _n\le x\pr - v \le x \pr |^2 \,dx + \int_{\Omega \setminus \Omega_\eta}|v_n \le x \pr - w \le x  \pr |^2\,dx \\
< \e_0^2/2 + \int _{\Omega \setminus \Omega_\eta} 4\le \alpha \le x \pr+ h\le x \pr \pr ^2 \,dx <\e_0^2.
\end{multline*}
Thus, contrary to \eqref{w-k_do_sprzecznosci}, $v_n \in \Psi\le t_0, u_0 \pr +
B_{L^2 \le \Omega , \R^N \pr} \le 0, \e_0 \pr$ for infinitely many $n\geq 1$.\hfill $\square$\\
\indent In order to get the weak tangency $(A_5)$ fix  $u\in K_0$ and define $G,H:\Omega \m \R^N, $ by
\[
G\le x \pr := F\le x, u \le x),\nabla u(x \pr \pr, \; H\le x \pr := T_C(u \le x\pr) \;\; \text{ for } x\in \Omega.
\]
The map $T_C(\cdot):C\m\R^N$ is lower semicontinuous (see \cite[Th. 4.2.2]{Aubin}), $G$ is measurable; hence $\Omega\ni x\m G(x)\cap H(x)\subset\R^N$ is measurable with nonempty values. By  the Kuratowski--Ryll-Nardzewski theorem, there is a measurable  $v:\Omega \to \R^N$ such that $v\le x \pr \in G \le x \pr \cap H \le x \pr$ for a.e. $x\in \Omega$. Clearly $v\in E$ and $v\in T_K(u)\cap F(t,u)$ since in view of \cite[Cor. 8.5.2]{Aubin} $T_K(u)=\{v\in E\mid v(x)\in T_C(u(x))\;\hbox{a.e. ion}\;\Omega\}.$\hfill $\square$

\begin{prop}\label{resolvent invariance I} For any $h>0$ the resolvent $J_h$ maps $K$ into $K_0$.\end{prop}
\noindent {\sc Proof:}
In view of \cite[Cor. 7.49]{Border}, $C$ is an intersection of countably many closed half-spaces containing it, i.e. $C=\bigcap_{n\geq 1}C_n$, where $C_n:=\left\{x\in \R^N\mid p_n\cdot x \leq a_n\right\}$ for some $p_n\in \R^N$ and $a_n\in \R$.
Thus, it is enough to show that $J_h(K_n)\subset K_n$, where
 $$K_n=\{u\in E\mid u(x)\in C_n\;\;\text{for a.a}\;\;x\in\Omega\}$$
 for every $n\geq 1$, since then
\[
J_h(K)=J_h\left(\bigcap_{n\geq 1}K_n\right)\subset \bigcap_{n\geq 1}J_h(K_n)\subset \bigcap_{n\geq 1}K_n = K\; (\footnote{Observe that in order to have that $\cap K_n\subset K$ one needs the {\em countable} collection of supporting functionals.}).
\]
But $J_h(K)\subset E_0$ so, eventually, $J_h(K)\subset K_0$.\\
\indent  Without loss of generality, we assume that $C=\left \{x\in \R^N\mid  p \cdot x \leq a\right \}$ for some $p\in\R^N $ and $a\in\R$ and $K=\{u\in E\mid p\cdot u(x)\leq a\;\;\text{for a.a}\;\;x\in\Omega\}$. Take $f\in K$ and put $u = J_h(f)$. By definition $u\in D(A)$ and
\[
u-hAu=f.
\]
Define $\bar{f}(x):= p\cdot f(x)$, $\bar{u}(x) := p\cdot u(x)$ for $x\in \Omega$.
Observe that $\bar{f}\leq a$ a.e.,  $\bar{u}\in H^2(\Omega)$  and for every $\xi \in H^1(\Omega)$
\[
\int_\Omega\bar{u}(x)\xi(x) \; dx + h\int_\Omega d(x)\nabla\bar{u}(x)\nabla\xi(x)\; dx = \int_\Omega\bar{f}(x)\xi(x)\; dx
\]
what yields
\[
\int_\Omega(\bar{u}(x)-a)\xi(x) \; dx= \int_\Omega(\bar{f}(x)-a)\xi(x)\; dx - h\int_\Omega d(x)\nabla(\bar{u}-a)(x)\nabla\xi(x)\; dx .
\]
Taking $\xi=(\bar{u}-a)_+:=\max\left\{0,\bar{u}-a\right\}$, we have $\xi\in H^1(\Omega)$  and $\nabla \xi = \chi\nabla (\bar{u}-a)$ by \cite[Cor. 1.3.6]{Cazenave}, where $\chi=\chi_{\left \{\bar{u}>a\right\}}$.
Therefore, for such $\xi$:
\[
0\leq \int_\Omega(\bar{u}-a)_+^2(x) \; dx=\int_\Omega(\bar{f}(x)-a)(\bar{u}- a)_+(x)\; dx -h\int_{\left\{\bar{u}>a\right\}} d(x)| \nabla(\bar{u}-a)_+(x)|^2\; dx \leq 0.
\]
As a result $\bar{u}\leq a$ a.e., that is $J_h(f)\in K$.\hfill $\square$

In view of Propositions \ref{properties of Phi} and \ref{resolvent invariance I} we get
\begin{thm}\label{main 3} If assumption $(D)$, $(F_1)$ -- $(F_3)$ are satisfied, the problem \eqref{problem II} has a solution.\hfill $\square$ \end{thm}

\subsection{The Neumann problem II}\label{Neumann II} We will establish the existence of steady state solutions to problem \eqref{problem 2bis}, i.e.
\be\label{problem IIbis}\begin{cases}-Lu+\Gamma u\in G(x,u)\;\;u(x)\in C,\;\; \text{a.e. on}\;\; \Omega,\\
\mbox{}\;\; \frac{\partial}{\partial {\bf n}}u=0\;\; \text{on}\;\; \partial\Omega,\end{cases}\ee
where $C\subset\R^N$ is compact and convex. We assume  $(D)$ and
\begin{enumerate}
\item[$(\Gamma)$] for all $i=1,...,N$, $\gamma_i=\gamma=(\gamma_1,...,\gamma_M)\in L^\infty(\Omega,\R^M)$; for any $i=1,...,N$;
    \item[$(G_1)$] $G:\bar\Omega\times C\multi\R^N$ is upper semicontinuous with compact convex values (\footnote{Note that $G$ is bounded, i.e. $\sup_{y\in G(x,u)}|y|<\infty$.});
        \item[$(G_2)$] $G$ is weakly tangent to $C$, i.e. $G(x,u)\cap T_C(u)\neq\emptyset$ for a.a. $x\in \Omega$ and all  $u\in C$.
        \end{enumerate}
Similarly as before we put
$$E=L^2(\Omega,\R^N),\;\; K:=\{u\in E\mid u(x)\in C\;\;\text{for a.a.}\;\;x\in\Omega\}.$$
Thus $(B_1)$ is satisfied. Let us define a continuous bilinear form
$$a(u,v)=\int_\Omega \left(d\,\nabla u\cdot \nabla v+\Gamma u\cdot v\right)\,dx,\;\; u,v\in H^1(\Omega,\R^N),$$
where $\nabla u \cdot\nabla v$ is the Frobenius product of derivatives (\footnote{Recall that if $A=[a_{ij}]$, $B=[b_{ij}]$ are $N\times M$ matrices, then $A\cdot B:=\sum_{i=1}^N\sum_{j=1}^Ma_{ij}b_{ij}$.}) and $\Gamma u\cdot v=\sum_{i=1}^N(\gamma\cdot\nabla u_i)v_i$. Observe that for any $v\in H^1(\Omega,\R^N)$ and $\eps>0$,
$$d_0\|\nabla v\|^2_{L^2}\leq a(v,v)-\int_\Omega\Gamma v\cdot v\,dx\leq a(v,v)+\|\gamma\|_{L^\infty}\left(\eps\|\nabla v\|^2_{L^2}+\frac{1}{4\eps}\|v\|^2_{L^2}\right),$$
where $d_0=\inf_{x\in\bar\Omega}|d(x)|$, in view of the so-called $\eps$-Cauchy inequality. Taking $0<\eps<d_0/2$ we get
\be\label{oszac} c\|\nabla v\|^2_{L^2}\leq a(v,v) +C\|v\|^2_{L^2}\ee
for some positive constants $c, C$. Therefore for all $v\in H^1(\Omega,\R^N)$
$$c\|v\|^2_{H^1}\leq a(v,v)+(c+C)\|v\|^2_{L^2}.$$
This implies that we are back in the situation of Remark \ref{ass comments} (1), (3) and putting
$$D(A)=\{u\in H^1(\Omega,\R^N)\mid \forall\,v\in H^1(\Omega,\R^N)\;\;\; a(u,v)=\la f,v\ra_{L^2}\;\; \text{for some}\;\;f\in L^2(\Omega,\R^N)\}$$the formula
$$Au=-f,\;\; u\in D(A),$$
where $f$ corresponds to $u$ as in the definition of $D(A)$, well defines a closed densely defined linear operator satisfying assumption $(B_2)$. Moreover $A$ is the generator of a $C_0$ semigroup $\{S(t)\}_{t\geq 0}$. The smoothness of the boundary $\partial\Omega$ and the standard regularity arguments imply that
$$D(A)=\left\{u\in H^2(\Omega,\R^N)\mid \left.\frac{\partial u}{\partial{\bf n}}\right|_{\partial\Omega}=0\right\}\;\;\;\text{and}\;\;\; Au:=Lu-\Gamma u,\;\,u\in D(A).$$
Now, for any $u\in K$, we put
$$\Phi(u):=\{v\in E\mid v(x)\in G(x,u(x))\;\; \text{for a.a.}\;\; x\in\Omega\}.$$
Following arguments from the proof of Proposition \ref{properties of Phi} we easily get  that $\Phi$ has properties $(B_3)$ and \eqref{Clarke tangency}.

In order to apply Theorem \ref{main 2} we need
\begin{prop}\label{resolvent invariance II} Condition  $(B_4)$ is satisfied.\end{prop}
\noindent {\sc Proof:} To this end we need use the $C_0$-semigroup structure. In view of Remark \ref{Euler 2} (2) we need to show that $K$ is semigroup invariant i.e. $S(t)u_0\in K$ for all $t\geq 0$ and $u_0\in K$. It is well known that
$$[0,+\infty)\ni t\mapsto u(t):=S(t)u_0$$ is the unique {\em mild} solution to the Cauchy initial value problem
\be\label{pp}\begin{cases} u'=Au,\;\; u\in E,\;\;t>0\\
u(0)=u_0.\end{cases}\ee
A function $v:[0,+\infty)\to E$ is a {\em strong} solution to \eqref{pp} if $v(0)=u_0$, $v(t)\in D(A)$ for $t\in (0,+\infty)$, $v\in W^{1,1}_{loc}((0,+\infty),E)$ and $v'(t)=Av(t)$ for a.a. $t\in (0,+\infty)$ ($v'(t)$ denotes the strong derivative of $v$ which exists a.a. since $v\in W^{1,1}_{loc}$). It is clear that each strong solution is a mild solution.\\
\indent Observe that $v$ is a strong solution to \eqref{pp} if and only if $w(t)\equiv e^{-\omega t}u(t)$, $t\in\I$, is a strong solution to $\dot w=(A-\omega I)w$. Hence, without loss of generality, we can assume that the form $a$ is nonnegative. In view of  \cite[Proposition II.2.5, Corollary III.2.4]{Show} (see also \cite[Th. 1.9.3]{Vrabie} and the classical result of Brezis in \cite[Th. 3.6]{Brezis}, comp. \cite[Prop. 5.1.1]{Cazenave}), if $f\in L^2([0,T], E)$, then there is a unique strong solution $v$ to the problem \eqref{pp} satisfying some additional conditions, which are irrelevant at the moment. As a consequence we see that $u$ is a strong solution to \eqref{pp}, i.e.
$$u'(t)=Au(t),\;\; t>0,\;\; u(0)=u_0\in K,$$
in the sense described above; moreover $u_0(x)\in C$ for a.a. $x\in\Omega$. Exactly as in the proof of Proposition \ref{resolvent invariance I}, using supporting functionals, we may assume that $C=\{x\in\R^N\mid p\cdot x\leq a\}$ for some $p\in\R^n$ and $a\in\R$. Put $\bar u:=p\cdot u$, i.e., $\bar u(t,x)=p\cdot u(t,x)$ and let $v:=\bar u-a$. Then $v(0)=p\cdot u_0-a$ and $v(t)\in H^1(\Omega)$ for $t>0$.  It is easy to see that $v:[0,+\infty)\to L^2(\Omega)$ is a strong solution to the following problem
$$v'(t)=\bar Av(t),$$where
$$\la \bar A\xi,\zeta\ra_{L^2(\Omega)}=-\bar a(\xi,\zeta):=-\int_\Omega d\,\nabla\xi\cdot \nabla\zeta\,dx-\int_\Omega(\gamma\cdot\nabla\xi)\zeta\,dx,\;\; \xi,\zeta\in H^1(\Omega).$$
Therefore $v(t)\in D(\bar A)$ for $t>0$ (as in the case of $A$, the domain of $\bar A$ consists of functions in $H^2(\Omega)$ whose normal outward derivative vanishes on $\partial\Omega$ in the sense of trace) and $v\in W^{1,1}_{loc}((0,+\infty),L^2(\Omega))$. For any $\xi\in H^1(\Omega)$ and $t>0$
$$\la v'(t),\xi\ra_{L^2}=-\bar a(v(t),\xi).$$
Now let us take $w(t)=v(t)_+$, i.e. $w(t):=\max\{\bar u(t)-a,0\}$, for $t\geq 0$. Then $w(t)\in H^1(\Omega)$, $t>0$, and by \cite[Prop. III.1.2]{Show},
$$\frac{1}{2}\frac{d}{dt}\|w(t)\|_{L^2}^2=\la w'(t),w(t)\ra_{L^2}^2=\la v'(t),w(t)\ra_{L^2}^2=-\bar a(v(t),w(t))=-\bar a(w(t),w(t)).$$
Using a counterpart of estimate \eqref{oszac} valid for $\bar a$ we see that for
$$\frac{1}{2}\frac{d}{dt}\|w(t)\|^2_{L^2}\leq -c\|\nabla w(t)\|^2_{L^2}+C\|w(t)\|_{L^2}^2\leq C\|w(t)\|^2_{L^2}.$$
By the Gronwall inequality we infer that $\|w(t)\|_{L^2}=0$ for all $t>0$, since $w(0)=0$ in $L^2(\Omega)$. It other words $u(t)\in K$ for all $t>0$.\hfill $\square$

\begin{thm}\label{main 4} If assumptions $(D)$, $(\Gamma)$, $(G_1)$ and $(G_2)$ are satisfied, then problem \eqref{problem IIbis} has a solution. \hfill $\square$
\end{thm}
\begin{rem} (1) Observe that theorems \ref{main 3} and \ref{main 4} are true if 0 belogs the constraint set $C$ and problems \eqref{problem II}, \eqref{problem IIbis} are subject to the Dirichlet condition. The only difference in proof is to see that if $0\in C$ and a functional $p$ supports $C$, i.e. $C\subset\{x\in\R^N\mid p\cdot x\leq a\}$, then $a\geq 0$. Hence $(p\cdot u-a)_+\in H^1_0(\Omega)$.\\
\indent (2) Let $C$ be convex compact and $Q=(0,l)^M \subset \R^M$ be an open cube. We will look for solutions to the problem \eqref{problem IIbis} with the periodic boundary conditions, i.e.
\begin{equation}
\begin{cases}
-Lu +\Gamma u \in G(x,u)\quad u(x)\in C \text{ a.e. on }Q,\\
u|_{\{x_i=0\}}=u|_{\{x_i=l\}},\quad \frac{\partial }{\partial x_i}\left|_{\{x_i=0\}}\right.=\frac{\partial}{\partial x_i}u\left|_{\{x_i =l\}}\right..
\end{cases}
\end{equation}
Assume (D), ($\Gamma$) (with $\Omega=Q$) and $d\in C^1_p(Q)$ (\footnote{This symbol stands for the restrictions to $Q$ of functions from $C^1(\R^M)$ which are $l$-periodic in each direction.}) and
\begin{enumerate}\label{prob.periodic}
\item[($P_1$)] $G:Q\times C \multimap \R^n$ is upper semicontinuous with compact convex values;
\item[($P_2$)] $G$ is weakly tangent to $C$, i.e. $G(x,u)\cap T_C(x)\neq\emptyset$ for all $x\in Q$ and $u\in C$.
\end{enumerate}
Let us put (\footnote{By  $H^k_p(Q)$ we denote the Sobolev space of $l$-periodic functions on the $M$-dimensional domain $Q$ ($k$ positive integer); see \cite[Chapter 5.10]{Robinson} for definitions and properties of $H^k_p(Q)$.}):
\begin{gather*}
E:=L^2(Q,\R^N),\quad K:=\{u\in E\mid u(x)\in C\text{ a.e. on }Q\},\\
D(A):=H^2_p(Q)^N\quad\text{and}\quad Au:=Lu -\Gamma u,\quad u\in D(A) .
\end{gather*}
We show, exactly as in section \ref{Neumann II}, that conditions $(B_1)$ -- $(B_4)$ are satisfied.  Thus, Theorem \ref{main 2} yields the existence of solutions to \eqref{prob.periodic}.
\end{rem}

\subsection{Some remarks to the Bernstein theory}
In a series of results in \cite{Granas} authors presented a modern approach to the so-called Bernstein theory for boundary value problems for  second order ordinary differential equations (see also \cite{Frigon},  \cite[II.7.4]{Dugundji-Granas}; for a numerous research afterwards see e.g. \cite{Tis} and bibliography therein). For the sake of completeness we formulate a model result \cite[Theorem 1.7]{Granas}.
\begin{thm}\label{bern} Suppose that $f:[0,T]\times\R\times\R\to\R$ is continuous such that
\begin{enumerate}
\item[{\em (i)}] (Sign condition) there is $R\geq 0$ such that $f(t,u,0)u\leq 0$ for $|u|>R$ and $t\in [0,T]$;
\item[{\em (ii)}] there are $a,b\geq 0$ so that $|f(t,u,v)|\leq av^2+b$ for $|u|\leq R$, $t\in [0,T]$ and $v\in\R$.
    \end{enumerate}
Then the problem
\be\label{ode 1}-u''=f(t,u,u')\ee subject to Dirichlet $(u(0)=u(T)=0)$, Neumann $(u'(0)=u'(T)=0)$ or periodic $(u(0)=u(T)$, $u'(0)=u'(T))$  boundary conditions has a solution in $C^2([0,T],\R)$ such that $|u(t)|\leq R$ for all $0\leq t\leq T$.\end{thm}

To illustrate our approach we will stay on the level of an ordinary differential inclusion and study the Dirichlet problem (the Neumann and periodic problems may be studied analogously) for
\be\label{ode 2}-u''\in \vp(t,u,u'),\ee
where
\begin{enumerate}
\item[$(\vp_1)$] $\vp:[0,T]\times\R^N\times\R^n\multi\R^N$ is upper continuous with compact convex values;
\item[$(\vp_2)$] there are $c, R\geq 0$ such that $\min_{y\in\vp(t,u,v)}y\cdot u\leq cR^2$ whenever $t\in [0,T]$, $\|u\|=R$, $v\in\R^N$ with $u\cdot v=0$;
    \item[$(\vp_3)$] $\vp$ is bounded on the strip $[0,T]\times C\times\R^N$, where $C:=\{y\in\R^N\mid |y|\leq R\}$.
\end{enumerate}
In order to apply Theorem \ref{main 1} let us put
$$\begin{gathered}E:=L^2((0,T),\R^N),\;\; K:=\{u\in E\mid |u(x)|\leq R\;\;\text{a.e. on}\;\; (0,T)\};\\
E_0:=C_0^1([0,T],\R^N)=\{u\in C^1([0,T],\R^N)\mid u(0)=u(T)=0\},\;\; K_0=K\cap E_0;\\
Au:=u''+cu\;\; \text{for}\;\; u\in D(A):=H^2((0,T),\R^N)\cap C_0([0,T],\R^N).\end{gathered}$$
Moreover define $\Phi:K_0\multi E$ by
$$\Phi(u)=\{v\in E\mid v(t)\in\vp(t,u(t),u'(t))-cu(t)\;\;\text{a.e. on}\;\; [0,T]\},\;\; u\in K_0.$$
Within this setting we see that conditions $(A_1)$, $(A_3)$ and $(A_4)$ are satisfied. As concerns $(A_2)$ note that for $u\in D(A)$ and $v\in V:=H_0^1((0,T),\R^N)$
$$\la Au,v\ra_{L^2}=-\int_0^Tu'\cdot v'\,dt+c\int_0^Tu\cdot v=-a(u,v),$$
where a bilinear form $a:V\times V\to\R$ given by
$$a(u,v)=\int_0^Tu'\cdot v'\,dt-c\int_0^Tu\cdot v\,dt,\;\; u, v\in V.$$
Thus for any $v\in V$
$$a(v,v)+c\|v\|^2_{L^2}=\|v\|^2_{H_0^1},$$
where $\|\cdot\|_{H_0^1}$ is the `short' norm in $V$. Hence we see that $A$ is the generator of the $C_0$ semigroup of linear operators on $E$ and conditions $(A_2)$, $(A_7)$ hold true, since the inclusion $D(A)\to E_0$ is compact. Condition $(A_6)$ may be shown as $(B_4)$ in Proposition \ref{resolvent invariance II}. Therefore we only need
\begin{prop} Condition $(A_5)$ is satisfied.\end{prop}
\noindent {\sc Proof:}  Let $u\in K_0$, i.e. $u\in C^1([0,T],\R^N)$, $u(0)=u(T)=0$, $|u(t)|\leq R$ for all $t\in [0,T]$. Let
$$X:=\{t\in [0,T]\mid |u(t)|=R\}.$$
It is clear that $X\subset (0,T)$. If $t\in X$, then $u(t)\cdot u'(t)=0$ since the function $(0,T)\ni s\mapsto |x(s)|^2$ takes maximum at $t$. Therefore there is $z\in\vp(t,u(t),u'(t))$ such that $z\cdot u(u)\leq cR^2=c|u(t)|^2$, i.e. $(z-cu(t))\cdot u(t)\leq 0$. Hence and by Remark \ref{examples} (2)
$$[\vp(t,u(t),u'(t))-cu(t)]\cap T_C(u(t))\neq\emptyset.$$
If $t\not\in X$, then $T_C(u(x))=\R^N$ and so
$$\forall\,t\in [0,T]\;\;\;[\vp(t,u(t),u'(t))-cu(t)]\cap T_C(u(t))\neq\emptyset.$$
Arguing as in the last part of the proof of Proposition \ref{properties of Phi} we produce a measurable $v\in\Phi(u)$ such that $v(t)\in T_C(u(t))$ for all $t\in [0,T]$. Then again by \cite[Cor. 8.5.2]{Aubin}, $v\in T_K(j(u))$.\hfill $\square$

According to Theorem \ref{main 1}, there is $u\in D(A)$ such that $-Au\in\Phi(u)$. This implies
\begin{thm}\label{problem IV} Under the above assumptions the problem \eqref{ode 2} has a strong solution $u\in K$.\hfill $\square$\end{thm}

\begin{cor}Suppose that $f:[0,T]\times\R^N\times\R^N\to\R^N$ is continuous such that
\begin{enumerate}
\item[{\em (i)}] there is $c,M\geq 0$ such that $u\cdot f(t,u,v)\leq cM^2$ for $t\in [0,T]$, $|u|=M$ and $v\in\R^N$ with $u\cdot v=0$ {\em (\footnote{Observe that if $N=1$, then this means that $u\cdot f(t,u,0)\leq cM^2$.})};
    \item[{\em (ii)}] $f$ is bounded on the strip  of the form $[0,T]\times D(0,M)\times\R^N$.
    \end{enumerate}
Then the problem
\be\label{ode 1}-u''=f(t,u,u')\ee subject to Dirichlet $(u(0)=u(T)=0)$, Neumann $(u'(0)=u'(T)=0)$ or periodic $(u(0)=u(T)$, $u'(0)=u'(T))$  boundary conditions has a solution in $C^2([0,T],\R)$ such that $|u(t)|\leq R$ for all $0\leq t\leq T$.\hfill $\square$
\end{cor}
The reader will easily formulate analogous results for elliptic PDE or partial differential inclusions. For instance one can get the generalization of the classical concerning the existence of steady states of the heat equation $u_t-\Delta u=g(u)$ subject to the Dirichlet boundary condition, where a continuous $g$ is such that for some positive $K,C$ one has $ug(u)\leq C|u|^2$ for $|u|\leq K$.

\subsection{Sub- and superharmonics; moving rectangles}
In this section we will discuss the Dirichlet problem
\be\label{moving}-\Delta u\in H(x,u,\nabla u),\;\;u|_{\partial\Omega}=0.\ee
Now we assume
\begin{enumerate}
\item[$(H_1)$] There are $\alpha,\beta\in C^1(\Omega,\R^N)\cap C(\bar\Omega,\R^N)$ such that $\alpha\leq\beta$  (\footnote{Here in below inequalities between vectors are understood in the componentwise sense.}), $\alpha$ is a weakly sub- and $\beta$ a weakly superharmonic, i.e. for all $\vp\in C^\infty_0(\Omega)$, $\vp\geq 0$,
    $$\int_\Omega\Delta\vp\,\alpha\,dx\geq 0,\;\; \alpha|_{\partial \Omega}\leq 0\;\; \text{and}\;\;\int_\Omega\Delta\vp\,\beta\leq 0,\;\; \beta|_{\partial\Omega}\geq 0;$$
\item[$(H_2)$] $H:C\times\R^{MN}\multi\R^N$ is a bounded upper semicontinuous map with convex compact values, where
    $$C:=\{(x,u)\in\Omega\times\R^N\mid \alpha(x)\leq u\leq\beta(x)\}.$$
    \end{enumerate}
Note that, for any $x\in\Omega$, the section $C(x):=\{u\in\R^n\mid (x,u)\in C\}$ is a cube in $\R^N$, hence the cone $T_{C(x)}(u)$ is determined in Remark \ref{examples} (1). The set $C$ maybe viewed as the graph of the {\em moving rectangles} $\Omega\ni x\mapsto C(x)$.

\indent For any $i=1,...,N$ let us introduce the lower and upper $i$-th `faces' of $C$:
$$C_i(x):=\{u\in\R^n\mid(x,u)\in C,\;u_i=\alpha_i(x)\},\;\; C^i(x):=\{u\in\R^N\mid (x,y)\in C,\;u_i=\beta_i(x)\}.$$
Further we assume that for all $i=1,...,N$ and $x\in\Omega$
\begin{enumerate}
\item[$(H_3)$] if $u\in C_i(x)$ and $v=(v_1,...,v_N)\in(\R^M)^N$ with $v_i=\nabla\alpha_i(x)$, then there is
    $y\in H(x,u,v)$ with $y_i\geq 0$;
\item[$(H_4)$] if $u\in C^i(x)$ and $v=(v_1,...,v_N)\in(\R^M)^N$ with $v_i=\nabla\beta_i(x)$, then there is
    $y\in H(x,u,v)$ with $y_i\leq 0$.
\end{enumerate}

Let
\begin{gather*} E:=L^p(\Omega,\R^N),\;\; E_0=C^1(\Omega,\R^N)\cap C_0(\bar\Omega,\R^N);\\
K:=\{u\in E\mid \alpha(x)\leq u(x)\leq \beta(x)\;\;\text{a.e. on}\;\;\Omega\},\;\; K_0=\{u\in E_0\mid \alpha(x)\leq \beta(x)\;\;\text{on}\;\;\Omega\};\\
D(A)=W_0^{1,p}(\Omega,\R^N)\cap W^{1,p}(\Omega,\R^N),\;\;Au=\Delta u,\;\; u\in D(A),\end{gather*}
where $p>M$.

\begin{thm} Under hypotheses $(H_1)$ -- $(H_4)$ problem \eqref{moving} has a strong solution in $K$.
\end{thm}
\noindent {\sc Proof:} It is clear that assumptions $(A_1)$ -- $(A_3)$ (with $\omega=0$) from section \ref{existence} are satisfied. Let us define $\Phi:K_0\multi E$ by
$$\Phi(u)=\{v\in E\mid v(x)\in H(x,u(x),\nabla u(x))\;\;\text{for}\;\;x\in\Omega\}.$$
As in Proposition \ref{properties of Phi} we check that assumption $(A_4)$ is also verified. Moreover $(A_7)$ holds true since $p>N$ and, thus, the inclusion $D(A)\to E_0$ is compact. We will check that $(A_5)$ and $(A_6)$ are true.\\
\indent Take $u\in K_0$ and let $X_i:=\{x\in\Omega\mid u(x)\in C_i(x)\}$, $X^i:=\{x\in\Omega\mid u(x)\in C^i(x)\}$, $i=1,..,N$. If $x\in X_i$ for some $i$, then $u_i(x)=\alpha_i(x)$ and $\nabla u_i(x)=\nabla\alpha_i(x)$ since $u_i-\alpha_i$ attains a minimum at $x$; similarly if $x\in X^i$, then $u_i=\beta_i(x)$ and $\nabla u_i(x)=\nabla\beta_i(x)$. Hence, by $(H_3)$ and $(H_4)$, if $x\in \bigcup_{i=1}^N(X_i\cup X_i)$, then $H(x,u(x),\nabla u(x))\cap T_{C(x)}(u(x))\neq \emptyset$. Otherwise, if $\alpha_i<u(x)<\beta_i(x)$ for all $i=1,..,N$, then $T_{C(x)}(u(x))=\R^N$. Takin into account
\cite[Cor. 8.5.2]{Aubin} we see that  $w\in T_K(j(u))$ if and only $w(x)\in T_{C(x)}(u(x))$ for a.a.  $x\in\Omega$. Arguing as in the last part of  the proof of Proposition \ref{properties of Phi} we see that $(A_5)$ is satisfied.\\
\indent Condition $(A_6)$ follows implicitly from \cite[Th. 16]{Kuiper}. Since we are in a special situation let us show a simple argument. Using $(H_1)$ and and the density arguments we see that if $v\in H_0^1(\Omega)$ and $v\geq 0$ a.e. on $\Omega$, then for any $i=1,..,N$,
$$\int_\Omega\nabla v\cdot\nabla \alpha_i\,dx\leq 0\;\; \text{and}\;\; \int_\Omega\nabla v\cdot\nabla\beta_i\geq 0.$$
Let $f\in K$, i.e. $\alpha\leq f\leq\beta$, and $u=J_h(f)$, where $h>0$. Then for any $v\in H_0^1(\Omega)$, $v\geq 0$,
$$\int_\Omega(u_i-\beta_i)v\,dx+h\int_\Omega(\nabla u-\nabla\beta_i)\cdot\nabla v\,dx=\int_\Omega(f_i-\beta_i)v\,dx-h\int_\Omega\nabla v\cdot\nabla\beta_i\,dx\leq 0.$$
Taking $v:=(u_i-\beta_i)_+$, we see that $v\in H_0^1(\Omega)$ since $\beta_i\geq 0$. Hence
$$\int_\Omega|v|^2\,dx\leq -\int_{\{u_i>\beta_i\}}|\nabla(u_i-\beta_i)|^2\,dx\leq 0$$
since $\nabla v=\chi_{\{u_i>\beta_i\}}\nabla(u_i-\beta_i)$. Thus $u_i\leq\beta$. Analogously we proof that $u_i\geq\alpha$ on $\Omega$ for all $i=1,...,N$. This means that $u\in K\cap D(A)\subset K_0$. Applying Theorem \ref{main 1} we end the proof.\hfill $\square$
\begin{rem} The existence of solutions to \eqref{moving} may be established by a direct use of arguments employed in the proof of Theorem \ref{main 1} since in this situation we can use some particular issues present in the problem. Given $i=1,...,N$ and $u\in W_0^{1,p}(\Omega,\R^N)$ let
$$\pi_i(u)=\begin{cases}\beta_i&\text{if}\;\; u_i<\beta_i,\\
u_i&\text{if}\;\; \alpha_i\leq u_i\leq\beta_i,\\
\alpha_i&\text{if}\;\; i_i<\alpha_i.\end{cases}$$
On can show that $\pi_i:{\mathcal E}:=W_0^{1,p}(\Omega,\R^N)\to W^{1,p}(\Omega)$ is well-defined and continuous. The map $\pi:=(\pi_1,...,\pi_N):{\mathcal E}\to{\mathcal E}$ is a retraction of
$\mathcal E$ onto $K\cap{\mathcal E}$. Note that $E_0\hookrightarrow {\mathcal E}$ and let
$$\Psi(u):=\Phi(\pi(u)),\;\;u\in E_0.$$ Taking into account that $0\in\rho(\Delta)$ we my consider the composition
$$\xi:E_0\stackrel{\Psi}{\multimap}E\stackrel{(-\Delta)^{-1}}{\longrightarrow}D(A)\to E_0.$$
This composition is a compact (at large: i.e. the range of $\xi$ is relatively compact) upper-semicontinuous map with compact convex values. By the Glicksberg-Fan theorem (the set-valued version of the Schauder fixed point principle) we gather that $\xi$ has a fixed point $u\in E_0$. Using $(H_3)$ and $(H_4)$ and the maximum principle one show that $u$ is located in $K$ and, therefore is a solution to \eqref{moving}.
\end{rem}


\end{document}